\documentclass[letterpaper, 10 pt, conference]{ieeeconf}  % Comment this line out if you need a4paper

\IEEEoverridecommandlockouts                              % This command is only needed if 
\usepackage{amsmath} % assumes amsmath package installed
\usepackage{amssymb}  % assumes amsmath package installed
\usepackage{mathtools}
\usepackage{url}
\usepackage{bbm}
\usepackage{courier}
\usepackage{algorithmicx}
\usepackage{subcaption}  % For subfigures
\usepackage{cleveref}
\usepackage{algorithm}
\usepackage[noend]{algpseudocode}

\usepackage{xcolor}
\usepackage{siunitx}
\usepackage{booktabs}
\usepackage[style=ieee, url=false,doi=false,eprint=false,natbib=false,mincitenames=1,maxcitenames=1,dashed=false,isbn=false, maxbibnames=99,backend=biber]{biblatex}
\usepackage{flushend}

\newcommand{\revboyd}[1]{{\color{black} #1}}

\def\baselinestretch{1.0}

\title{\LARGE \bf
Resource Allocation under Stochastic Demands\\
using Shrinking Horizon Optimization
}

\author{
  Alexandros E.~Tzikas,$^{1}$ Nazim Kemal Ure,$^{1}$ Mansur Arief,$^{1}$\\
  Mykel J. Kochenderfer,$^{1}$ and Stephen P. Boyd$^{2}$
  \thanks{$^{1}$A. E. Tzikas (corresponding author), N. Kemal Ure, M. 
  Arief, and M. J. Kochenderfer are with the Department of Aeronautics and Astronautics, Stanford University, Stanford, CA 94305, U.S.A.
        {\tt\small \{alextzik, ure, ariefm, mykel\}@stanford.edu}}%
    \thanks{$^{2}$S. P. Boyd is with the Department of Electrical Engineering, Stanford University, Stanford, CA 94305, U.S.A.
        {\tt\small boyd@stanford.edu}}%
}

\begin{document}

\maketitle
\thispagestyle{empty}
\pagestyle{empty}

%%%%%%%%%%%%%%%%%%%%%%%%%%%%%%%%%%%%%%%%%%%%%%%%%%%%%%%%%%%%%%%%%%%%%%%%%%%%%%%%
\begin{abstract}
    We consider the problem of optimally allocating a limited number of resources across time to maximize revenue under stochastic demands. This formulation is relevant in various areas of control, such as supply chain, ticket revenue maximization, healthcare operations, and energy allocation in power grids. We propose a bisection method to solve the static optimization problem and extend our approach to a shrinking horizon algorithm for the sequential problem. The shrinking horizon algorithm computes future allocations after updating the distribution of future demands by conditioning on the observed values of demand. \revboyd{We illustrate the method on a simple synthetic example with jointly log-normal demands, showing that it achieves performance close to a bound obtained by solving the prescient problem. The code for all experiments can be found at: \url{https://github.com/sisl/sequential_resource_allocation}. }
\end{abstract}

% \begin{keywords}
% planning under uncertainty, probability and statistical methods, optimal control, distribution steering.
% \end{keywords}

%%%%%%%%%%%%%%%%%%%%%%%%%%%%%%%%%%%%%%%%%%%%%%%%%%%%%%%%%%%%%%%%%%%%%%%%%%%%%%%%
\section{Introduction}

% Stephen's feedback:
% \begin{itemize}
%     \item Mention newsvendor problem
%     \item Mention applications of problem formulation.
%     \item Make the distinction between the static problem and the sequential MPC-like policy. Our contribution is mostly on the latter where we will use log normal distributions in experiments. The former could be a contribution depending on what Kemal can find in the literature.
% \end{itemize}
The problem of optimally allocating a finite quantity of resources over a discrete time horizon to satisfy stochastic demands is a fundamental challenge in control theory and operations research. This formulation can be viewed as a multi-period generalization of the classic newsvendor problem~\cite{qin2011newsvendor}, extending the single-period inventory decision to a sequence of interdependent choices. The practical applications of this model are extensive and span areas such as supply chain management~\cite{choi2008mean}, dynamic pricing with fixed capacity~\cite{monahan2004dynamic}, operational planning in healthcare~\cite{bavafa2019managing}, and the dispatch of energy resources in power grids~\cite{marchi2019multi}. In each case, the objective is to maximize the total revenue, subject to a constraint on the total available resources.

In this paper, we first analyze the static multi-period resource allocation problem. In this formulation, the allocation decisions for the entire horizon are determined at the beginning of time. This approach yields a convex optimization problem, for which we can find an efficient solution using a bisection method on the dual variable. However, the solution to the static problem is inherently limited; it relies solely on the marginal distributions of demand for each time period and fails to incorporate temporal correlations or the information revealed as demands are sequentially realized. It is an open-loop policy that \revboyd{does not} adapt to observed outcomes.

To overcome these limitations, we propose a sequential, adaptive policy based on a shrinking horizon optimization framework, similar to model predictive control (MPC)~\cite{schwenzer2021review}. Under this policy, at every time, an allocation decision is made only for the current time. At each subsequent time, the optimization problem is \revboyd{re-solved} for the remaining (shrinking) time horizon. Crucially, the joint distribution of future demands is updated by conditioning on the history of observed demands. This closed-loop approach allows the allocation strategy to dynamically adapt to new information as it becomes available. While this sequential method is a heuristic and does not guarantee global optimality, it is designed to systematically improve upon the static solution by leveraging the observed data.

The primary contribution of this work is the formulation and empirical evaluation of this shrinking horizon policy. We demonstrate its superiority against the static, open-loop solution. The evaluation is conducted on synthetic data where demands are drawn from a jointly log-normal distribution, a model chosen to reflect the correlated and non-negative nature of demands often found in practice. The results show that the adaptive method achieves higher revenue by effectively exploiting the information revealed over time. \revboyd{Its performance is close to a bound obtained by solving the prescient problem. }

The remainder of this paper is organized as follows. In \cref{sec:related_work} we discuss prior work. In \cref{sec:static} we describe the static optimization problem to maximize expected revenue and in \cref{sec:static_sol} we show how to solve it. In \cref{sec:sequential} we describe the sequential problem of maximizing revenue and in \cref{sec:sequential_sol} we demonstrate how to obtain a policy for it. Finally, in \cref{sec:sims} we describe our simulation set-up and results. We conclude the paper in \cref{sec:conclusion}.

\section{Related Work}\label{sec:related_work}

The static multi-period resource allocation problem is an extension of the classic newsvendor model~\cite{petruzzi1999pricing}. We consider a portfolio of decisions constrained by a total resource limit. This class of problems has been extensively studied in network revenue management~\cite{talluri2004revenue}, but not under log-normal stochastic models. While foundational, the classic newsvendor model's assumption of a known demand distribution has been challenged by modern, data-rich environments. Recent advancements have thus focused on the data-driven newsvendor problem, where decisions are made directly from demand samples without assuming a specific parametric distribution~\cite{levi2015data}, \cite{see2010robust}. This line of work bridges the gap between classic stochastic models and modern adaptive methods, though often still focusing on a single-period setting~\cite{caro2010inventory}. Even when extended to multi-period settings, these static, open-loop policies are inherently suboptimal as they fail to adapt to the sequential arrival of information~\cite{den2015dynamic}. 

Adaptive policies have been proposed for the multi-period resource allocation problem. We mention some works that we believe capture the general trends of prior research. Kim et al. propose a policy assuming the uncertainty is expressed as a discrete set of possible scenarios. Contrary to our formulation, the random quantities do not appear in the optimization objective \cite{kim2015optimal} and a discrete set of scenarios is used. Wang et al. propose an approximate dynamic programming approach under a complicated, partially observable information pattern \cite{wang2010optimal}. Dynamic programming has been used to solve the problem under the assumption of independent and identically distributed demands \cite{khang2000optimality}. Scenario trees, whose complexity scales exponentially with the number of nodes, have been proposed to solve discretized versions of the problem \cite{densing2013dispatch}. 
Approaches based on mixtures of experts have also been suggested \cite{wang2023solving}. Our method is different from the aforementioned works in the following ways: we propose and solve an abstraction of the sequential allocation problem that is widely applicable. Our problem formulation is highly descriptive. Our algorithm is fast, scales well with the problem horizon, and is applicable under any model for the demands.

% Addressing the sequential nature of the problem optimally requires solving a dynamic program (DP). In this context, the state would track remaining resources and the current time to derive a state-dependent allocation policy. Unfortunately, this approach is thwarted by the curse of dimensionality, as the state space becomes intractably large~\cite{puterman2014markov}. To circumvent this, various approximation methods have been proposed. Approximate dynamic programming (ADP) and, more recently, deep reinforcement learning (RL) have emerged as powerful paradigms, replacing the true value function with a learned approximation~\cite{bertsekas2019reinforcement}. While successful in complex domains, these methods can require extensive offline training and their performance often depends on sophisticated model architecture and hyperparameter tuning~\cite{agarwal2021deep}. 

Our shrinking horizon approach is a form of MPC, a well-known technique for approximating dynamic programming solutions. While MPC has been applied to supply chain problems, many implementations rely on certainty-equivalent formulations that use point forecasts~\cite{sarimveis2008dynamic}. Our method differs by directly incorporating the stochastic nature of the problem at each step, updating the full joint distribution of future demands based on past observations. This allows our policy to explicitly react to temporal demand correlations, a significant challenge that often complicates other adaptive approaches.

\section{The Static Problem}\label{sec:static}

We consider a sequence of random variables $d_1, \dots, d_T \in \mathbb{R}$, where $T$ is the discrete-time horizon of the problem and the random variables represent the unknown demands at every time. A demand $d_t$ means that $d_t$ units or items were demanded at time $t$. Because demands are nonnegative, it always holds that
\begin{equation}
    d_1, \dots, d_T \geq 0.
\end{equation}
We do not assume that the demands $d_1, \dots, d_T$ are independent or identically distributed. In many applications, demands at different times are correlated and their marginal distributions can be significantly different. 

We assume that we have knowledge of the joint distribution of $d_1, \dots, d_T$ and that for any time $t=1, \dots, T$, we can compute the \revboyd{marginal} cumulative distribution function (CDF) $F_t(x)=:\mathbb{P}\left( d_t \leq x \right)$ for any $x\in \mathbb{R}$. \revboyd{Equivalently, we assume that we can compute the quantile functions $F_t^{-1}(y)$ for $y \in \left[0, 1 \right]$ and $t=1, \dots, T$}. This is possible, for example, when we can sample $(d_1, \dots, d_T)$ from their joint distribution or when the joint distribution is known in closed form. 

The static resource allocation problem is 
\begin{equation} \label{eq:main}
\begin{aligned}
    \max_{a_1, \dots, a_T} \qquad &\mathbb{E}_{d_1, \dots, d_T} \left[ \sum_{t=1}^T p_t \min \left(d_t, a_t \right)\right]\\
    \mathrm{s.t.}\qquad &\sum_{t=1}^T a_t = L\\
    &a_1, \dots, a_T \geq 0
\end{aligned},
\end{equation}
where the variables $a_1, \dots, a_T \in \mathbb{R}$ are the allocated units or items for each time and $p_t > 0$ is the price per unit at time $t$. We are maximizing the expected revenue by optimally allocating the resources across time. The first constraint implies that there is a total of $L >0$ units to be allocated across time and the second constraint does not allow for a negative allocation on any given time. 

Problem \eqref{eq:main} is a convex program. To see this, we proceed as follows. The objective of \eqref{eq:main} can be rewritten as
\begin{equation}\label{eq:objective}
    \sum_{t=1}^T p_t \mathbb{E}_{d_t} \left[ \min \left( d_t, a_t \right)\right],
\end{equation}
where only the marginal distributions of the random variables, $d_1, \dots, d_T$ are used.
For every realization $\Tilde{d}$, the function $\min \left( \Tilde{d}, a_t \right)$ is concave in $a_t$. The expectation can be viewed as a sum with nonnegative weights of terms $\min \left( \Tilde{d}, a_t \right)$ for different values $\Tilde{d}$. Therefore, $\mathbb{E}_{d_t} \left[ \min \left( d_t, a_t \right)\right]$ is concave in $a_t$. Because $p_t>0$, we finally get that the function in \eqref{eq:objective} is concave in $(a_1, \dots, a_T)$.

\section{The Static Solution}\label{sec:static_sol}

We show how to obtain the solution to problem \eqref{eq:main}. First, we find an equivalent expression for the objective. Then, we leverage duality to solve the problem.

\subsection{An Equivalent Form for the Objective of \eqref{eq:main}}

Consider the term $\mathbb{E}_{d_t} \left[ \min \left(d_t, a_t \right) \right]$. We can write
\begin{equation}
    \begin{aligned}
        \mathbb{E}_{d_t} \left[ \min \left(d_t, a_t \right) \right] &= \int_0^\infty \mathbb{P}\left( \min\left( d_t, a_t\right) \geq x \right)dx\\
        &=\int_0^{a_t} \mathbb{P}\left( \min\left( d_t, a_t\right) \geq x \right)dx\\
        &=\int_0^{a_t} \mathbb{P}\left( d_t \geq x \right) dx.
    \end{aligned}
\end{equation}
The first equality is a well-known result in probability \cite{billingsley2017probability}. The second equality follows because $\min\left(d_t, a_t\right) \leq a_t$. Finally, the third equality follows because for any $0\leq x\leq a_t$
\begin{equation}
    \min\left( d_t, a_t\right) \geq x \Leftrightarrow  d_t \geq x.
\end{equation}

By defining the survival functions
\begin{equation}
S_t(x) =: \mathbb{P}\left( d_t > x\right) = 1 - \mathbb{P}\left(d_t \leq x \right),
\end{equation}
we can write the objective of \eqref{eq:main} as
\begin{equation}
    \sum_{t=1}^T p_t \int_0^{a_t} S_t(x)dx,
\end{equation}
where the optimization variables $a_t$ appear on the top limit of the integrals.

\subsection{Optimality Conditions}
The Lagrangian of problem \eqref{eq:main} is given by
\begin{equation}
\begin{aligned}
    L(a_1, \dots, a_T, \lambda_1, \dots, \lambda_T, \nu) =\\ \sum_{t=1}^T p_t \int_0^{a_t} \mathbb{P}\left( d_t \geq x\right)dx 
    + \sum_{t=1}^T \lambda_t a_t - \nu \left( \sum_{t=1}^T a_t - L \right),
\end{aligned}
\end{equation}
where $\nu$ is the dual variable for the resource limit constraint and, for every $t$, $\lambda_t \geq 0$ is the dual variable for the non-negativity constraint $a_t \geq 0$.

The Karush--Kuhn--Tucker (KKT) optimality conditions \cite{boyd2009convex} dictate that the primal-dual optimal point, $(a_1^\star, \dots, a_T^\star, \lambda_1^\star, \dots, \lambda_T^\star, \nu^*)$, satisfies
\begin{subequations}\label{eq:KKT}
    \begin{align}
        p_t \dfrac{\partial}{\partial a_t} \left(\int_0^{a_t} \mathbb{P}\left( d_t \geq x\right) dx)\right)\vert_{a_t^\star} + \lambda_t^\star &= \nu^\star,\label{eq:KKT_1}\\
        \sum_{t=1}^T a_t^\star &= L,\label{eq:KKT_2}\\
        \lambda_t^\star a_t^\star &= 0\label{eq:KKT_3},\\
        \lambda_t^\star \geq 0,\quad a_t^\star &\geq 0.
    \end{align}
\end{subequations}
Using the Leibniz integral rule \cite{border2016differentiating}, the first condition \eqref{eq:KKT_1} can be written as 
\begin{equation}\label{eq:inter}
    p_t \mathbb{P}\left( d_t \geq a_t^\star \right) + \lambda_t^\star =\nu^\star.
\end{equation}
Eq.~\eqref{eq:inter} offers an intuitive explanation for the solution $(a_1^\star, \dots, a_T^\star)$. For times $t$ that the price per item, $p_t$, is high, the allocation $a_t^\star$ is large. To see this, simply ignore $\lambda_t^\star$. However, the relationship between price $p_t$ and allocation $a_t^\star$ is not linear. Rather, it depends on the distribution of $d_t$.

We can obtain an expression for $a_t^\star$ as a function of $\nu^\star$. We suppose that $\mathbb{P}\left( d_t \leq \delta \right) > 0$ for any $\delta >0$. 
\begin{itemize}
    \item Suppose $a_t^\star >0$. Then, by \eqref{eq:KKT_3}, $\lambda_t^\star =0$ and, by \eqref{eq:inter}, we get that $\mathbb{P}\left( d_t \geq a_t^\star \right)= \nu^\star / p_t < 1$, by the assumption above.
    \item Suppose $a_t^\star =0$. Then $\mathbb{P}\left( d_t \geq a_t^\star \right)=1$ and, since $\lambda_t^\star \geq 0$, by \eqref{eq:inter}, we get that $\nu^\star / p_t \geq \mathbb{P}\left( d_t \geq 0 \right)$, i.e., $\nu^\star / p_t \geq 1$.
\end{itemize}
Therefore, we obtain that
\begin{equation}\label{eq:a_t}
    a_t^\star = \begin{cases}
        F_t^{-1}(1-\dfrac{\nu^\star}{p_t}),\quad &\dfrac{\nu^\star}{p_t} < 1,\\
        0,\quad &\mathrm{otherwise}
    \end{cases},
\end{equation}
where we can obtain $\nu^\star$ through \eqref{eq:KKT_2}. Obviously by \eqref{eq:inter}, it holds that $\nu^\star \geq 0$.

Notice from \eqref{eq:a_t} that $\sum_{t=1}^T a_t^\star$ is a decreasing function of $\nu^\star$. 
Furthermore, suppose that for every $t$, it holds that $\nu^\star / p_t \geq 1$. Then $a_t^\star = 0$ for all $t$, which is not possible because $L>0$. Therefore, there exists a $t$ such that $\nu^\star / p_t < 1$, which implies 
\begin{equation}
    \nu^\star \in \left[0, \max \{p_1, \dots, p_T \}\right].
\end{equation}
We can solve the system \eqref{eq:KKT_2}-\eqref{eq:a_t} using bisection on $\nu^\star$. We can use $0$ and $\max \{p_1, \dots, p_T \}$ as the initial lower and upper limits for $\nu^\star$, respectively. This procedure is outlined in \Cref{alg:bisection}. For an accuracy of $\epsilon$ for $\nu^\star$, we require $\mathcal{O} \left( \log \left( \max\{p_1, \dots, p_T \}/\epsilon \right) \right)$ iterations. The only information required about the distribution of the demands are the quantile functions $F_t^{-1}(\cdot)$. 
These are straight-forward to compute for various distributions, such as Gaussian or Laplace. Even for Gaussian mixture models, this is easy using a root-finding algorithm \cite{kochenderfer2019algorithms}. We call the allocations produced by \Cref{alg:bisection} the \textit{static} solution.

In some applications, the allocations need to be integer-valued. We do not impose such a constraint in problem \eqref{eq:main}. However, after obtaining the solution $(a_1^\star, \dots, a_T^\star)$ using \Cref{alg:bisection}, we can round the allocation at each time to the closest integer and slightly modify them to ensure they sum to $L$. \revboyd{We can also easily do a local optimization around the rounded values. We note that with these slight modifications the duality gap is still bounded by a small number and not by a factor growing with $T$ because there is only one equality constraint \cite{udell2016bounding, udell2013maximizing}.}

It is important to note that we can solve problem \eqref{eq:main} using \Cref{alg:bisection} at the beginning of time and then allocate the optimal number of units $a_1^\star, \dots, a_T^\star$ over time. However, this approach can be sub-optimal because it only uses the marginal distributions of $d_1, \dots, d_T$ and ignores the correlation between the demands across time. Furthermore, the demands $d_1, \dots, d_T$ are realized sequentially, which \Cref{alg:bisection} does not leverage to improve performance.

\begin{algorithm}[t]
\caption{Solution of the static resource allocation problem \eqref{eq:main} using bisection}\label{alg:bisection}
\begin{algorithmic}[1]
\State \textbf{Inputs:} horizon $T$, prices per item $p_1, \dots, p_T >0$, resource allocation limit $L>0$, CDFs $F_t(u) = \mathbb{P}\left( d_t \leq u\right)$ for all $t=1, \dots, T$
\State \textbf{Parameters:} accuracy $\epsilon$ for solution $\nu^\star$ 
\State $N_\mathrm{iterations} \leftarrow \lceil \log \left( \max \{p_1, \dots, p_T \} / \epsilon \right) \rceil$ 
\State $\nu_\mathrm{low} \leftarrow 0$
\State $\nu_\mathrm{up} \leftarrow \max \{p_1, \dots, p_T \}$
\For{$k = 1, \dots, N_\mathrm{iterations}$}
    \State Obtain the midpoint $\nu_\mathrm{mid} \leftarrow \dfrac{\nu_\mathrm{low}+\nu_\mathrm{up}}{2}$
    \State Obtain $\hat{a}_1, \dots, \hat{a}_T$ through \eqref{eq:a_t}
    \If{$\sum_{t=1}^T \hat{a}_t \geq L$}
        \State $\nu_\mathrm{low} \leftarrow \nu_\mathrm{mid}$
    \Else
        \State $\nu_\mathrm{up} \leftarrow \nu_\mathrm{mid}$
    \EndIf
\EndFor
\State \Return optimal allocations $\hat{a}_1, \dots, \hat{a}_T$, optimal dual variable $\nu_\mathrm{mid}$
\end{algorithmic}
\end{algorithm}

\section{The Sequential Problem}\label{sec:sequential}
We consider the case where the allocations need not all be determined at the beginning of time. At any time $t$, we assume that we have observed the demands $d_1, \dots, d_{t-1}$ and already made the allocations $a_1, \dots, a_{t-1}$, which cannot be altered. Our goal at time $t$ is to determine the allocation $a_t$ using the available information, which includes the realized demands, $d_1, \dots, d_{t-1}$, the allocations made so far, $a_1, \dots, a_{t-1}$, and the joint distribution of the demands $d_1, \dots, d_T$.
\revboyd{We assume that, for any time $t=1, \dots, T$ we have access to the marginal CDFs and quantile functions for the demands, conditioned on the previously observed values, i.e., $1, \dots, t-1$ for every $t$.}

\section{The Sequential Policy}\label{sec:sequential_sol}

We propose a policy based on shrinking horizon optimization for the sequential resource allocation problem. Our policy takes into account the realized demands so far to update the distribution of future demands through conditioning. More precisely, at time $0$, we solve problem \eqref{eq:main} to obtain the solution $(a_{1\mid 0}^\star, \dots, a_{T\mid 0}^\star)$, where the subscript denotes that the solution does not depend on any observed values of demand, since none have been observed yet.

After assigning the allocation $a_{1\mid 0}^\star$ for time $1$, we observe the realized value for the demand $d_1$. To compute the future allocations, we now solve
\begin{equation} \label{eq:rec_1}
\begin{aligned}
    \max_{a_2, \dots, a_T} \qquad &\mathbb{E}_{d_2, \dots, d_T\mid d_1} \left[ \sum_{t=2}^T p_t \min \left(d_t, a_t \right)\right]\\
    \mathrm{s.t.}\qquad &\sum_{t=2}^T a_t = L - a_{1\mid 0}^\star\\
    &a_2, \dots, a_T \geq 0
\end{aligned},
\end{equation}
which has the solution $(a_{2\mid 1}^\star, \dots, a_{T\mid 1}^\star)$. In the objective of \eqref{eq:rec_1} we have conditioned on the observed value of the demand.
We assign the allocation $a_{2\mid 1}^\star$ at time $2$ and then observe the realized demand $d_2$. To make our approach clear, we mention that at time $3$, we will solve
\begin{equation} \label{eq:rec_2}
\begin{aligned}
    \max_{a_3, \dots, a_T} \qquad &\mathbb{E}_{d_3, \dots, d_T\mid d_1, d_2} \left[ \sum_{t=3}^T p_t \min \left(d_t, a_t \right)\right]\\
    \mathrm{s.t.}\qquad &\sum_{t=3}^T a_t = L - a_{1\mid 0}^\star - a_{2\mid1}^\star\\
    &a_3, \dots, a_T \geq 0
\end{aligned}.
\end{equation}

Our procedure is outlined in \Cref{alg:receding}. \Cref{alg:receding} requires knowing the conditional quantile functions for the demands, \revboyd{conditioned on the previously observed values}. As in \cref{sec:static_sol}, we can round the allocation at each time to the closest integer, if this is required by the application. Our approach is a sophisticated heuristic and need not be optimal. We call the allocations produced by \Cref{alg:receding}, the \textit{sequential} solution.

\revboyd{As a sidenote, in the series of problems \eqref{eq:rec_1}-\eqref{eq:rec_2}, we use the conditional mean in the objective. We could have also used the maximum a posteriori (MAP) value. Our policy is a heuristic and which option works best depends on the application.}

\subsection{The Case of Independent Demands}\label{sec:indep}
In the case of independent demands $d_1, \dots, d_T$, conditioning on realized demands does not alter the distribution of future demands. For example, the marginal distribution of  
$d_2, \dots, d_T$ is the same as the conditional distribution $d_2, \dots, d_T \mid d_1$. This implies that 
\begin{equation}\label{eq:sols_eq}
    (a_{2\mid 0}^\star, \dots, a_{T\mid 0}^\star) = (a_{2\mid 1}^\star, \dots, a_{T\mid 1}^\star).
\end{equation}
To see why \eqref{eq:sols_eq} is true, let
\begin{equation}
\begin{aligned}
    R_{2:} &=: \mathbb{E}_{d_2, \dots, d_T}\left[ \sum_{t=2}^T p_t \min\left( d_t, a_{t\mid 0}^\star\right) \right]\\
    &=\mathbb{E}_{d_2, \dots, d_T \mid d_1}\left[ \sum_{t=2}^T p_t \min\left( d_t, a_{t\mid 0}^\star\right) \right].
\end{aligned}
\end{equation}
The equality follows from independence. We suppose, for the sake of contradiction, that there exists an elementwise non-negative $(a_2, \dots, a_T)$ such that 
\begin{subequations}
\begin{align}
    \mathbb{E}_{d_2, \dots, d_T}\left[ \sum_{t=2}^T p_t \min\left( d_t, a_t\right) \right] > R_{2:},\\
    \sum_{t=2}^T a_t = L - a_{1\mid 0}^\star.\label{eq:sum_constr}
\end{align}
\end{subequations}
Then, it holds that $(a_{1\mid 0}^\star, a_2, \dots, a_T)$ is a feasible point for problem \eqref{eq:main} and obtains a larger objective value than $(a_{1\mid 0}^\star, \dots, a_{T\mid 0}^\star)$. However, this is a contradiction. Therefore, for any elementwise non-negative  $(a_2, \dots, a_T)$ that satisfies \eqref{eq:sum_constr} it holds that
\begin{equation}
\begin{aligned}
    \mathbb{E}_{d_2, \dots, d_T \mid d_1}\left[ \sum_{t=2}^T p_t \min\left( d_t, a_t\right) \right] \leq \\
    \mathbb{E}_{d_2, \dots, d_T \mid d_1}\left[ \sum_{t=2}^T p_t \min\left( d_t, a_{t\mid 0}^\star\right) \right]
\end{aligned}
\end{equation}
Therefore, \eqref{eq:sols_eq} holds, assuming a unique optimal point. Using recursion, we can prove that for any $\tau=1, \dots, T-1$:
\begin{equation}
    (a_{\tau+1\mid 0}^\star, \dots, a_{T\mid 0}^\star)  = (a_{\tau+1\mid \tau}^\star, \dots, a_{T\mid \tau}^\star).
\end{equation}

\revboyd{Thus, when the demands are independent, the \textit{sequential} solution obtained by our shrinking horizon policy reduces to the \textit{static} solution.}

\begin{algorithm}[t!]
\caption{Shrinking horizon procedure for the sequential resource allocation problem \eqref{eq:main}}\label{alg:receding}
\begin{algorithmic}[1]
\State \textbf{Inputs:} horizon $T$, prices per item $p_1, \dots, p_T >0 $, resource allocation limit $L>0$
\State Solve problem \eqref{eq:main} using \Cref{alg:bisection} to obtain $(a_{1\mid 0}^\star, \dots, a_{T\mid 0}^\star)$
\State Assign the allocation $a_{1\mid 0}^\star$ at time $1$ and observe the realized demand $d_1$
\For{$\tau=2, \dots, T$}
    \State Solve problem 
    \begin{equation}
\begin{aligned}
    \max_{a_\tau, \dots, a_T} \qquad &\mathbb{E}_{d_\tau, \dots, d_T\mid d_1, \dots, d_{\tau-1}} \left[ \sum_{t=\tau}^T p_t \min \left(d_t, a_t \right)\right]\\
    \mathrm{s.t.}\qquad &\sum_{t=\tau}^T a_t = L - \sum_{t=1}^{\tau-1} a_{t\mid t-1}^\star\\
    &a_\tau, \dots, a_T \geq 0
\end{aligned}
\end{equation}
using \Cref{alg:bisection} to obtain $(a_{\tau \mid \tau-1}^\star, \dots, a_{T \mid \tau-1}^\star)$.
\State Assign the allocation $a_{\tau \mid \tau-1}^\star$ at time $\tau$ and observe the realized demand $d_\tau$
\EndFor
\end{algorithmic}
\end{algorithm}
\section{Simulations}\label{sec:sims}
We compare the performance of the \textit{static} sand \textit{sequential} solutions on synthetic data. First, we describe the model used for the demands and the generation process for the data. Second, we describe an upper bound for the obtained revenue and a simple baseline. Third, we show the equivalence of the \textit{static} and \textit{sequential} solutions in the case of independent demands. Finally, we demonstrate the superiority of the \textit{sequential} solution in the case of non-independent demands. 

\subsection{Model for the Demands}
We model the demands $d_1, \dots, d_T$ as jointly log-normal random variables. Equivalently, we model 
$\log d_1, \dots, \log d_T$ as a jointly Gaussian random vector with mean $\mu$ and covariance $\Sigma$, i.e., 
\begin{equation}
    \left( \log d_1, \dots, \log d_T \right) \sim \mathcal{N}(\mu, \Sigma).
\end{equation}
The log-normal model is common in operations research, supply chain management, and financial models \cite{horvath1959applicability, dufresne2004log}. It is particularly suitable to model demands, because demands are non-negative, empirically exhibit right-skewness, and are correlated across time. The correlations across time can be incorporated elegantly, because the log demands are modeled as jointly Gaussian.

We now describe how to compute probabilities of demand under the log-normal model. We note that because $\mathbb{P}\left( d_t \leq x\right) = \mathbb{P}\left( \log d_t \leq \log x \right)$,
\begin{equation}
    F_t(x) = \Phi \left( \dfrac{\log x-\mu_t}{\sqrt{\Sigma_{t,t}}}\right),
\end{equation}
where $\mu_i$ is the $i$-th entry of vector $\mu$, $\Sigma_{i,j}$ is the $(i,j)$-th element of matrix $\Sigma$, and $\Phi(\cdot)$ is the CDF of the standard normal distribution. The conditional distribution of $d_t$ given $d_{t-1}, \dots, d_1$ can be written as
\begin{equation}\label{eq:cond}
\begin{aligned}
    \mathbb{P}\left( d_t \leq x \mid d_{t-1} = x_{t-1}, \dots, d_1 = x_1 \right) =\\
    \mathbb{P} \left( \log d_t \leq \log x \mid \log d_{t-1} = \log x_{t-1}, \dots, \log d_1 =\log x_1\right),
\end{aligned}
\end{equation}
because $\log$ is a monotone transformation. Assuming the conditional mean of $\log d_t$ is $\mu_{t\mid t-1}$ and its conditional standard deviation is $\sigma_{t \mid t-1}$, we get that
\begin{equation}
\begin{aligned}
    \mathbb{P}\left( d_t \leq x \mid d_{t-1} = x_{t-1}, \dots, d_1 = x_1 \right) =\\
    \Phi\left( \dfrac{\log x - \mu_{t\mid t-1}}{\sigma_{t \mid t-1}}\right).
\end{aligned}
\end{equation}
We can compute $\mu_{t\mid t-1}$ and $\sigma_{t \mid t-1}$ explicitly using conditioning on jointly Gaussian distributions \cite[Section 3.9]{kochenderfer2022algorithms}. 

Although modeling the log demands as jointly Gaussian is adequate for many applications, we can model them using a Gaussian mixture model, which is a universal approximator \cite{goodfellow2016deep}. \revboyd{In the case of the Gaussian mixture model, we can again efficiently evaluate the conditional quantile functions.}

\subsection{The Synthetic Data}
We describe the process of generating the synthetic data. We select prices $p_t$ between $10$ and $100$. We select the means of the demands $d_t$ between $20$ and $100$, while the standard deviations of the demands are between $10$ and $30$. The first moments of the demands are mapped appropriately to the means and standard deviations of the log demands using standard formulas \cite{balakrishnan2016continuous}. 

The correlations between the log demands are selected differently in our experiments depending on whether the demands are independent or not. 
For the case of independent demands we set the correlations to be $0$ for log demands of different times.
For the case of non-independent demands, we set the correlations between $\log d_t$ and $\log d_\tau$ to be (roughly) between $-0.7$ and $0.7$ for $\tau \neq t$. To do so, we first create a symmetric matrix with diagonal entries equal to $1$ and off-diagonal entries in $\{ -0.7, 0.7\}$. We then project to the set of positive semi-definite matrices to obtain the correlation matrix.

We set the limit $L$ to be between $30\%$ and $60\%$ of the sum of mean demands. The synthetic data is generated once and remains fixed across the simulation trials for each set of experiments below.

\subsection{A Revenue Upper-Bound}
Suppose we know the realized sequence of demands, $\Tilde{d_1}, \dots, \Tilde{d_T}$. Then, we can calculate the optimal allocations for this realization by solving the convex problem
\begin{equation} \label{eq:uppberbound}
\begin{aligned}
    \max_{a_1, \dots, a_T} \qquad &\sum_{t=1}^T p_t \min \left(\Tilde{d}_t, a_t \right)\\
    \mathrm{s.t.}\qquad &\sum_{t=1}^T a_t = L\\
    &a_1, \dots, a_T \geq 0
\end{aligned}.
\end{equation}
We call the solution of \eqref{eq:uppberbound}, the \textit{oracle} allocation. This allocation cannot be computed in a causal manner and is only known after observing the realized demands. Nevertheless, it serves as a useful upper-bound of performance. 

\subsection{The Roll-Forward Baseline}
\revboyd{
We describe a simple causal baseline. At time $1$, we set the allocation to be $L/T$. At each subsequent time $t$, we set the allocation to be the last observed demand, i.e., $a_t = d_{t-1}$, if the remaining allocation is larger then $d_{t-1}$. Otherwise, we simply set $a_t$ to be the remaining allocation. We call this baseline, the \textit{roll-forward} allocation, because it just rolls the last observed demand forward as the next allocation.}

\subsection{Independent Demands}
\begin{figure*}[tb]
\centering

% Row 1: Independent
\begin{subfigure}[t]{0.4\textwidth}
\centering
\includegraphics[width=\linewidth]{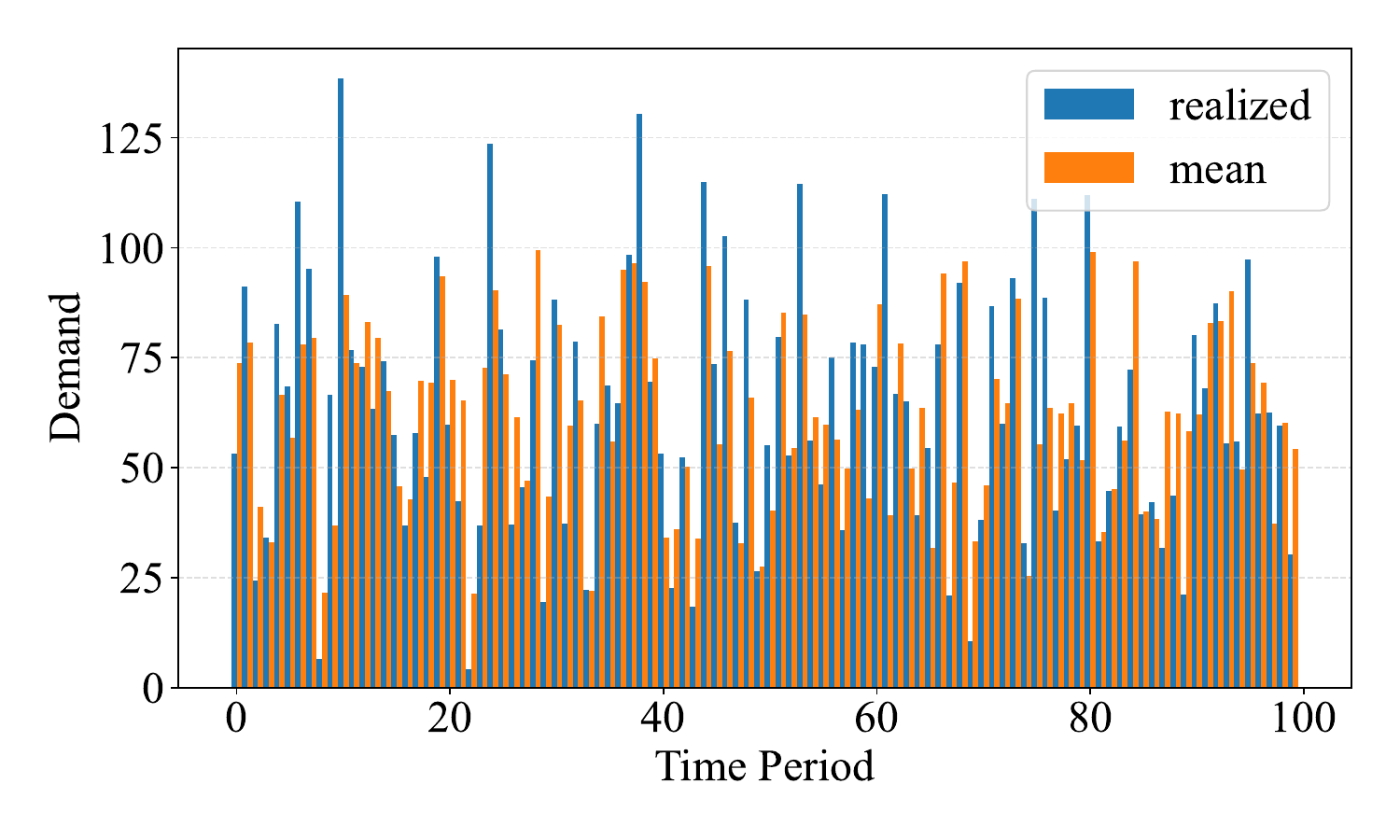}
\caption{Realized and expected demands across time.}
\label{fig:demand_indep}
\end{subfigure}
\hfill
\begin{subfigure}[t]{0.4\textwidth}
\centering
\includegraphics[width=\linewidth]{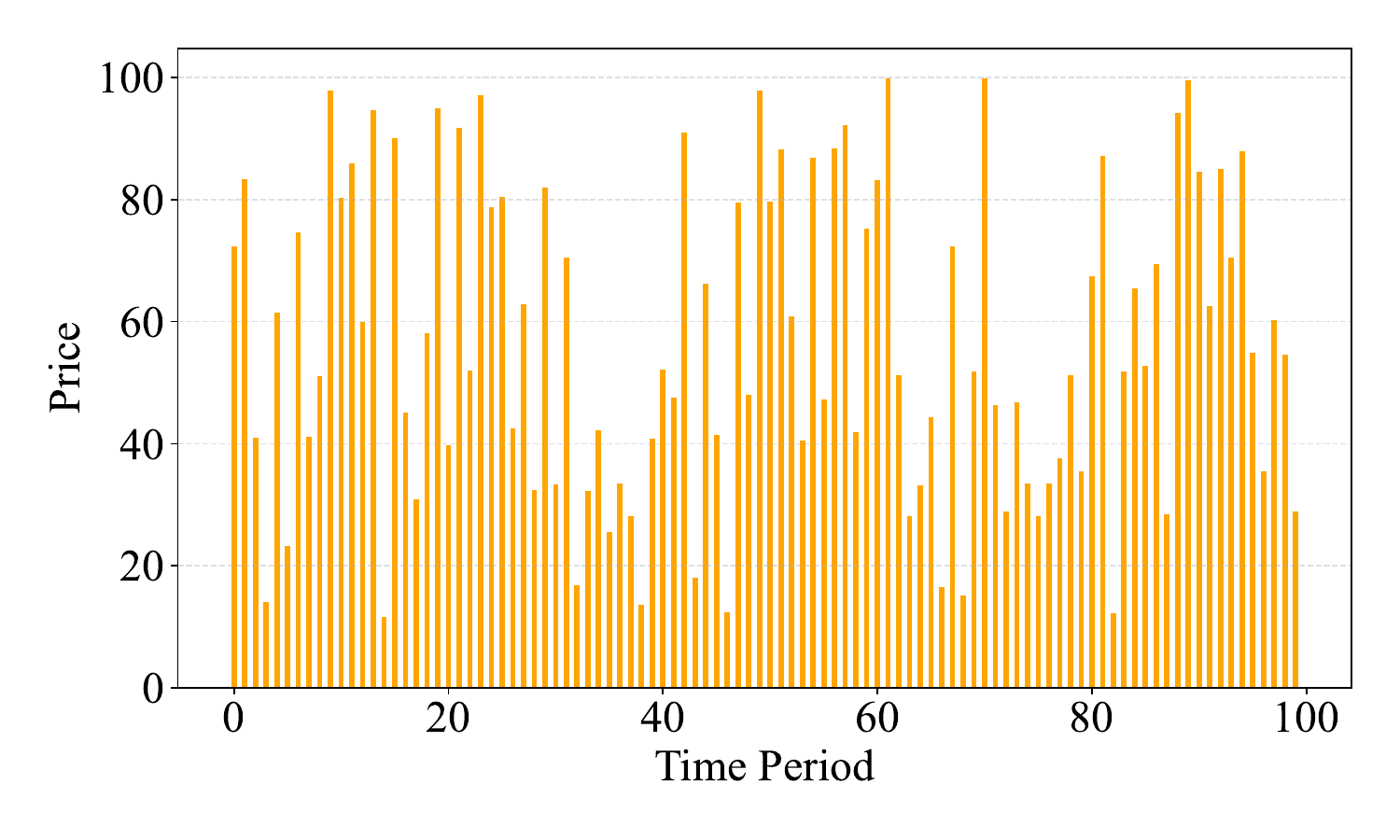}
\caption{Prices (per unit) across time.}
\label{fig:prices_indep}
\end{subfigure}

\vspace{0.45em}

% Row 2: Correlated
\begin{subfigure}[t]{0.4\textwidth}
\centering
\includegraphics[width=\linewidth]{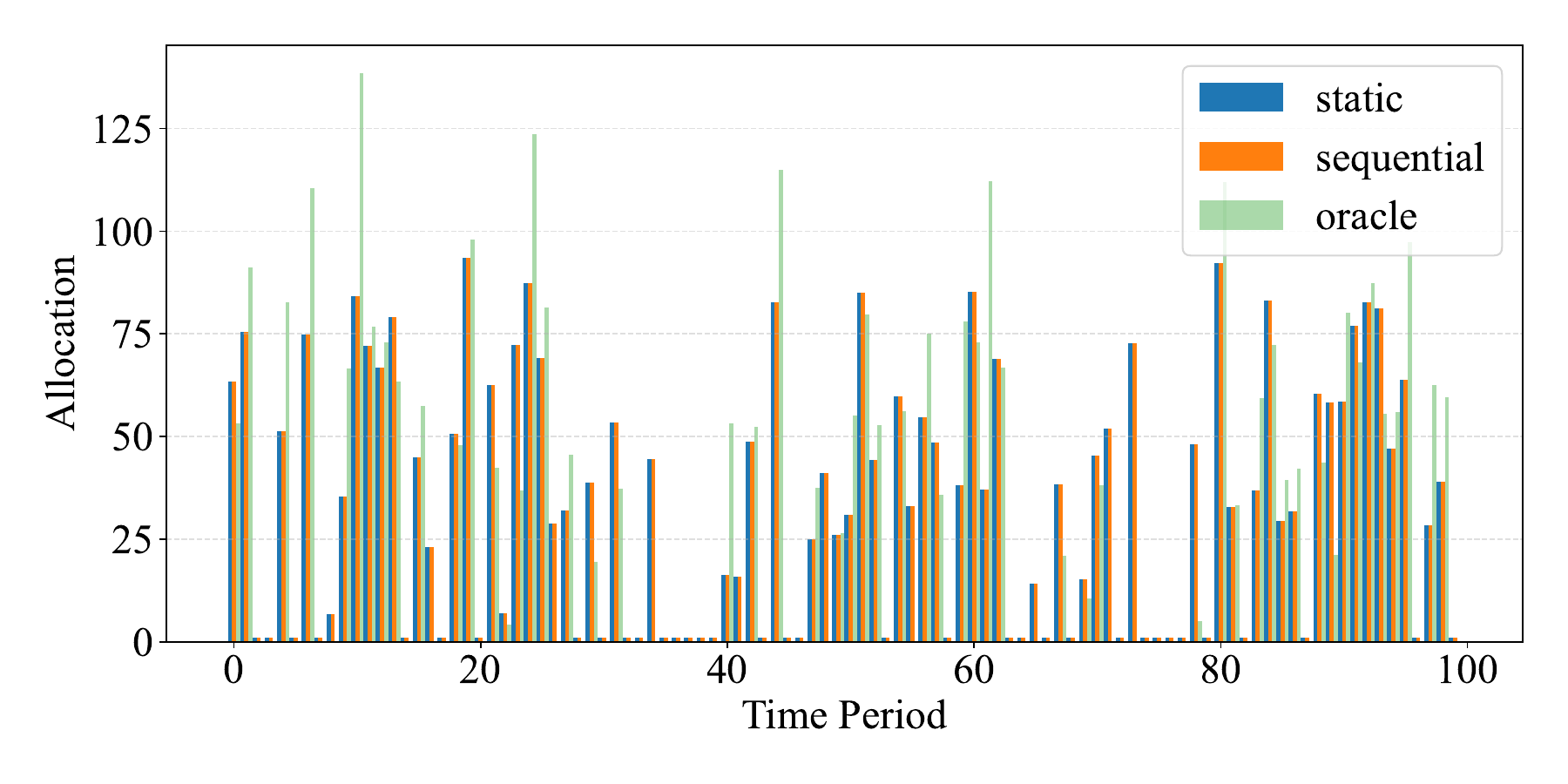}
\caption{The \textit{oracle}, \textit{sequential}, and \textit{static} allocations across time. The \textit{static} and \textit{sequential} allocations are the same for the case of independent demands.}
\label{fig:alloc_indep}
\end{subfigure}
\hfill
\begin{subfigure}[t]{0.4\textwidth}
\centering
\includegraphics[width=\linewidth]{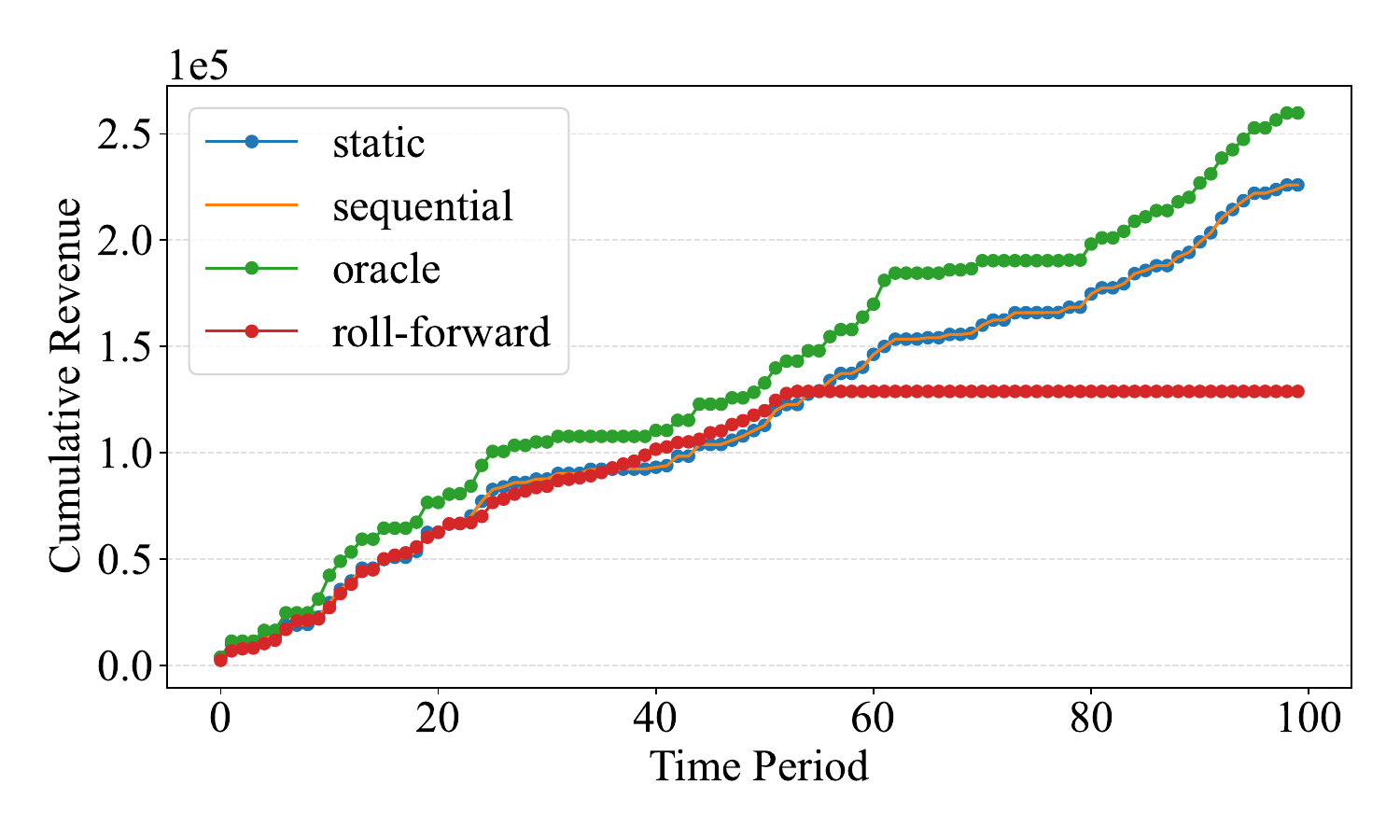}
\caption{The cumulative revenue over time for the \textit{oracle}, \textit{sequential}, \textit{static}, and \textit{roll-forward}  allocations. The cumulative revenue is the same for the \textit{static} and \textit{sequential} allocations.}
\label{fig:cum_ret_indep}
\end{subfigure}

\caption{Results for the case of independent demands.}
\label{fig:indep}
\end{figure*}

We set the correlation between $\log d_t$ and $\log d_\tau$ to be $0$ for any $\tau \neq t$. This implies that the log demands are independent. As functions of the log demands, the demands are also independent in this case. By \cref{sec:indep}, we know that the \textit{static} and \textit{sequential} solutions must match. This is verified in our simulation, as shown in \Cref{fig:alloc_indep}, for a particular realization of demands, which is shown in \Cref{fig:demand_indep}. In this case, the two approaches obtain the same revenue, which is smaller than the revenue obtained by the \textit{oracle} allocation, as shown in \Cref{fig:cum_ret_indep}. However, the \textit{oracle} allocation cannot be computed in a causal manner. \revboyd{The \textit{roll-forward} method, which is causal, performs the worst, because it does not adequately capture the inter-dependencies between the demands. It also allocates the total $L$ in far fewer time periods, which explains the constant cumulative revenue in later time periods. In \Cref{fig:alloc_indep}, we do not plot the \textit{roll-forward} allocation because this can be easily inferred by the realized demands in \Cref{fig:demand_indep}.}

\subsection{Non-Independent Demands}

We set the correlations between $\log d_t$ and $\log d_\tau$ to be (roughly) between $-0.7$ and $0.7$ for $\tau \neq t$. We first look at the performance of the methods for a single trial, i.e., a single realized sequence of demands, and then provide aggregate results over $100$ trials for different time horizons $T$.

\begin{figure*}
\centering

% Row 1: Independent
\begin{subfigure}[t]{0.4\textwidth}
\centering
\includegraphics[width=\linewidth]{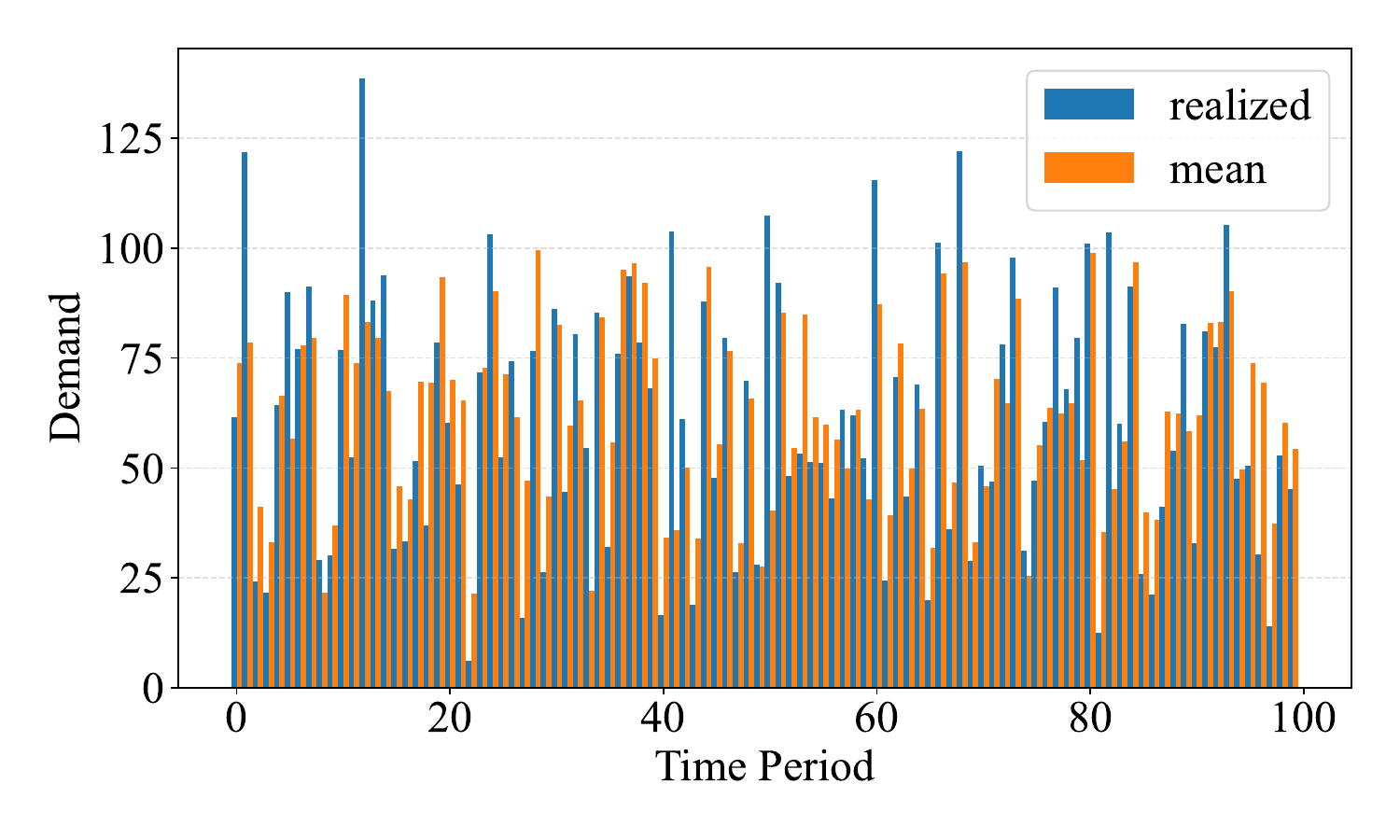}
\caption{Realized and expected demands across time.}
\label{fig:demand_dep}
\end{subfigure}
\hfill
\begin{subfigure}[t]{0.4\textwidth}
\centering
\includegraphics[width=\linewidth]{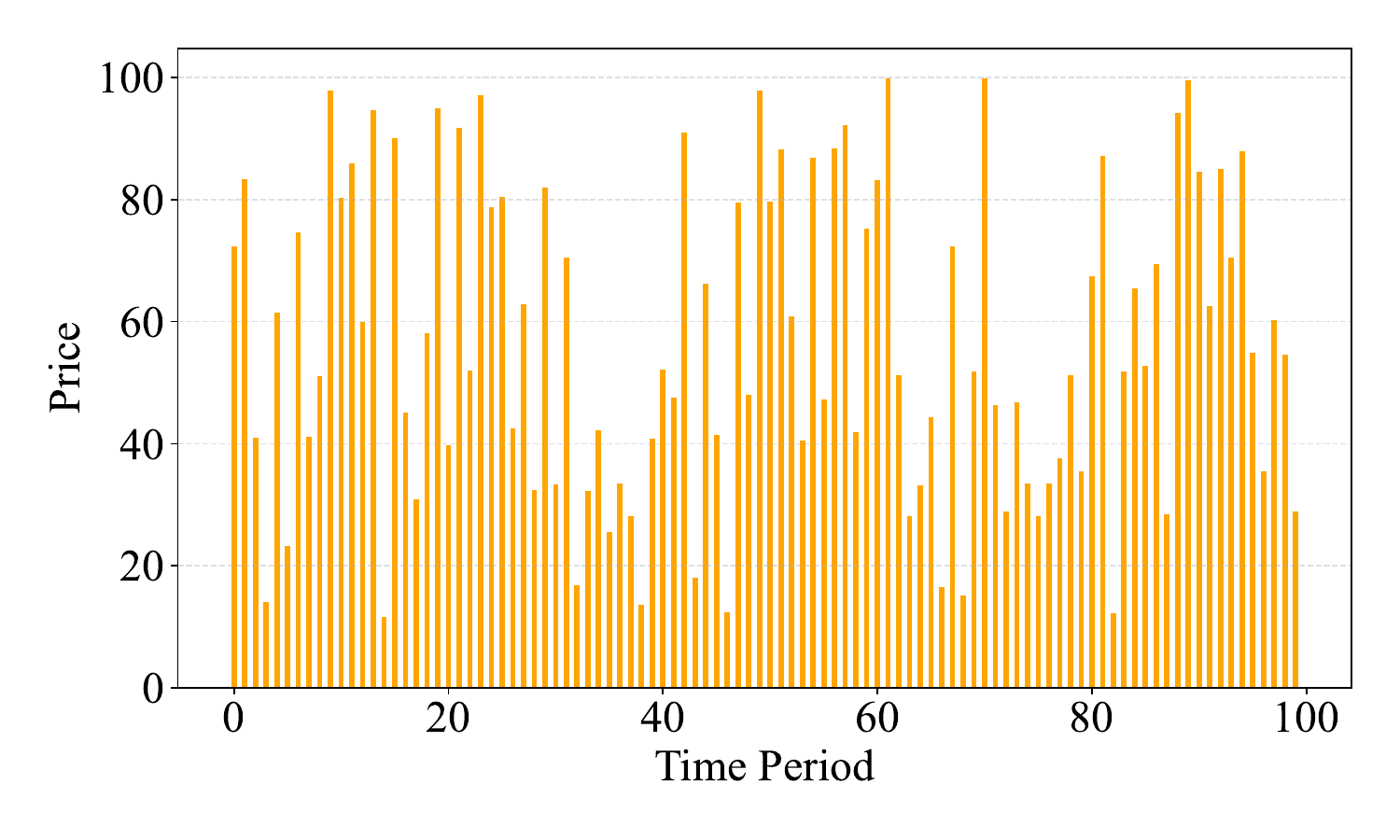}
\caption{Prices (per unit) over time.}
\label{fig:prices_dep}
\end{subfigure}

\vspace{0.45em}

% Row 2: Correlated
\begin{subfigure}[t]{0.4\textwidth}
\centering
\includegraphics[width=\linewidth]{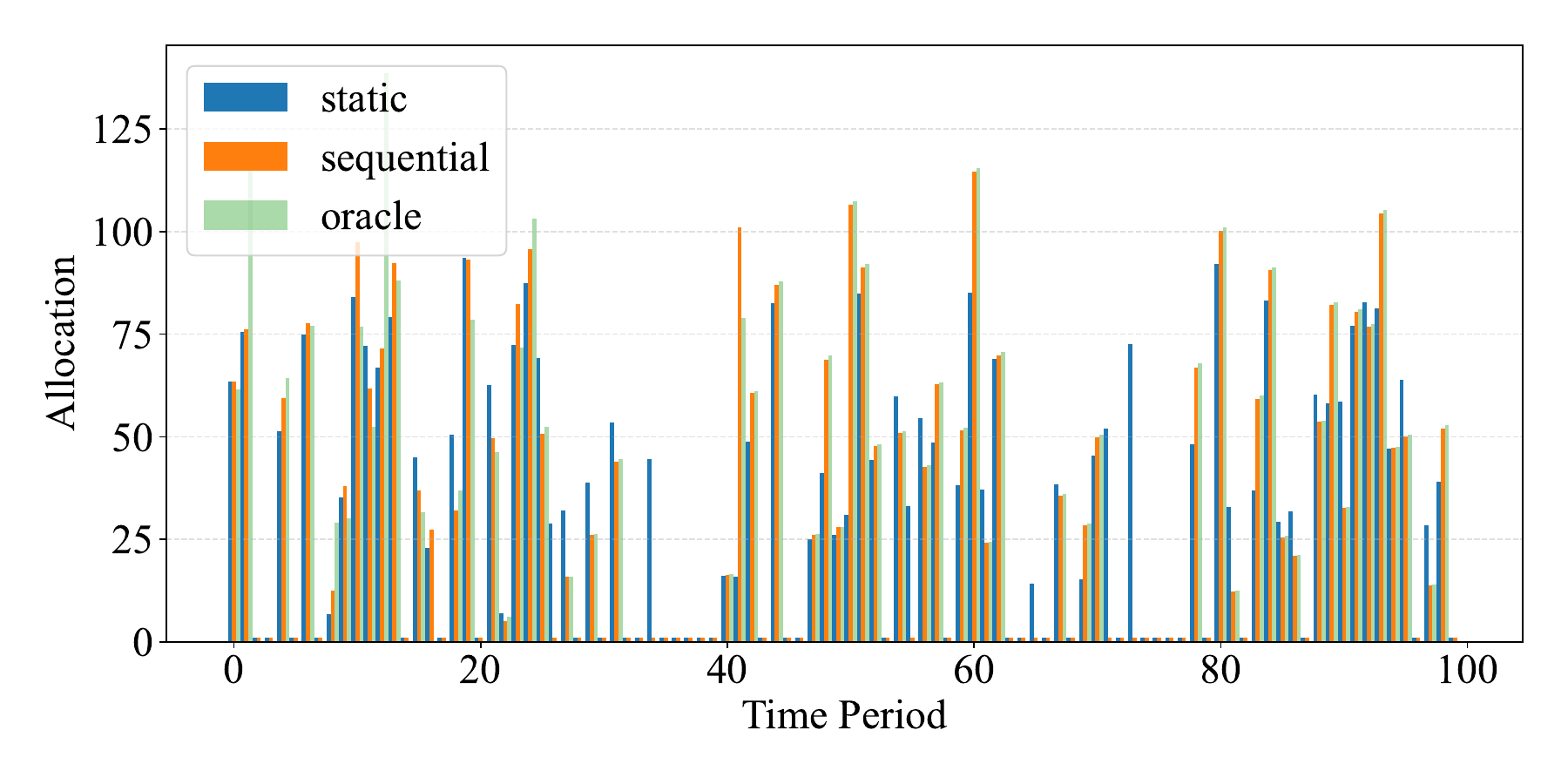}
\caption{The \textit{oracle}, \textit{sequential}, and \textit{static} allocations across time. The \textit{sequential} allocation is close to the \textit{oracle} allocation.}
\label{fig:alloc_dep}
\end{subfigure}
\hfill
\begin{subfigure}[t]{0.4\textwidth}
\centering
\includegraphics[width=\linewidth]{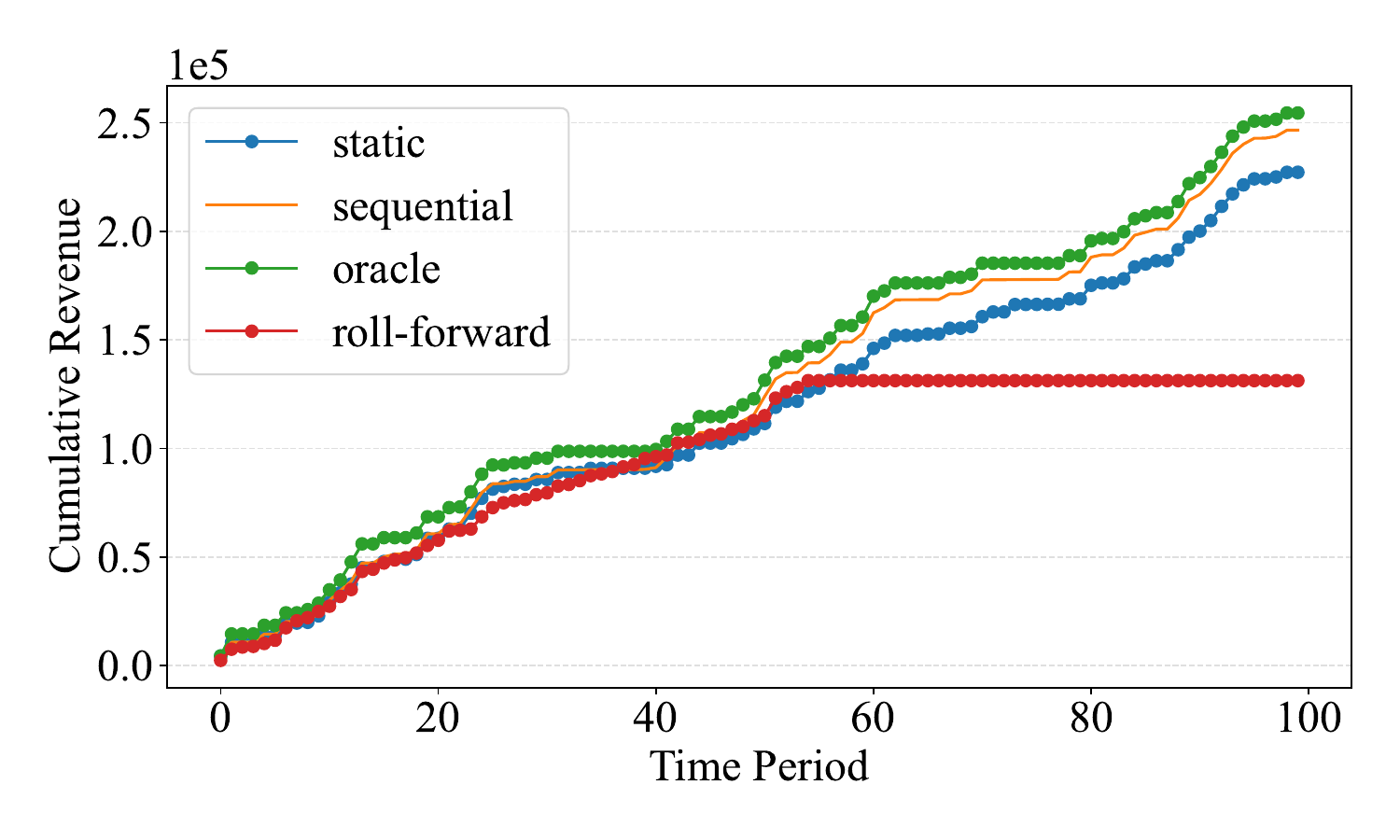}
\caption{The cumulative revenue over time for the \textit{oracle}, \textit{sequential}, \textit{static}, and \textit{roll-forward}  allocations. }
\label{fig:cum_ret_dep}
\end{subfigure}

\caption{Results for the case of non-independent demands.}
\label{fig:dep}
\end{figure*}

Our results for a single trial and $T=100$ are included in \Cref{fig:dep}. We observe that the \textit{sequential} allocation is closer to the oracle allocation than the \textit{static} one. As a result, the cumulative revenue of the \textit{sequential} allocation is higher than the cumulative revenue of the \textit{static} one. We note that we expect the gap between the revenue of the \textit{sequential}
and \textit{static} allocations to shrink as the correlations between the demands get closer to $0$, i.e., when the demands are close to being independent. 

\begin{figure*}
  \centering

  \begin{subfigure}[t]{0.4\linewidth}
    \centering
    \includegraphics[width=\linewidth]{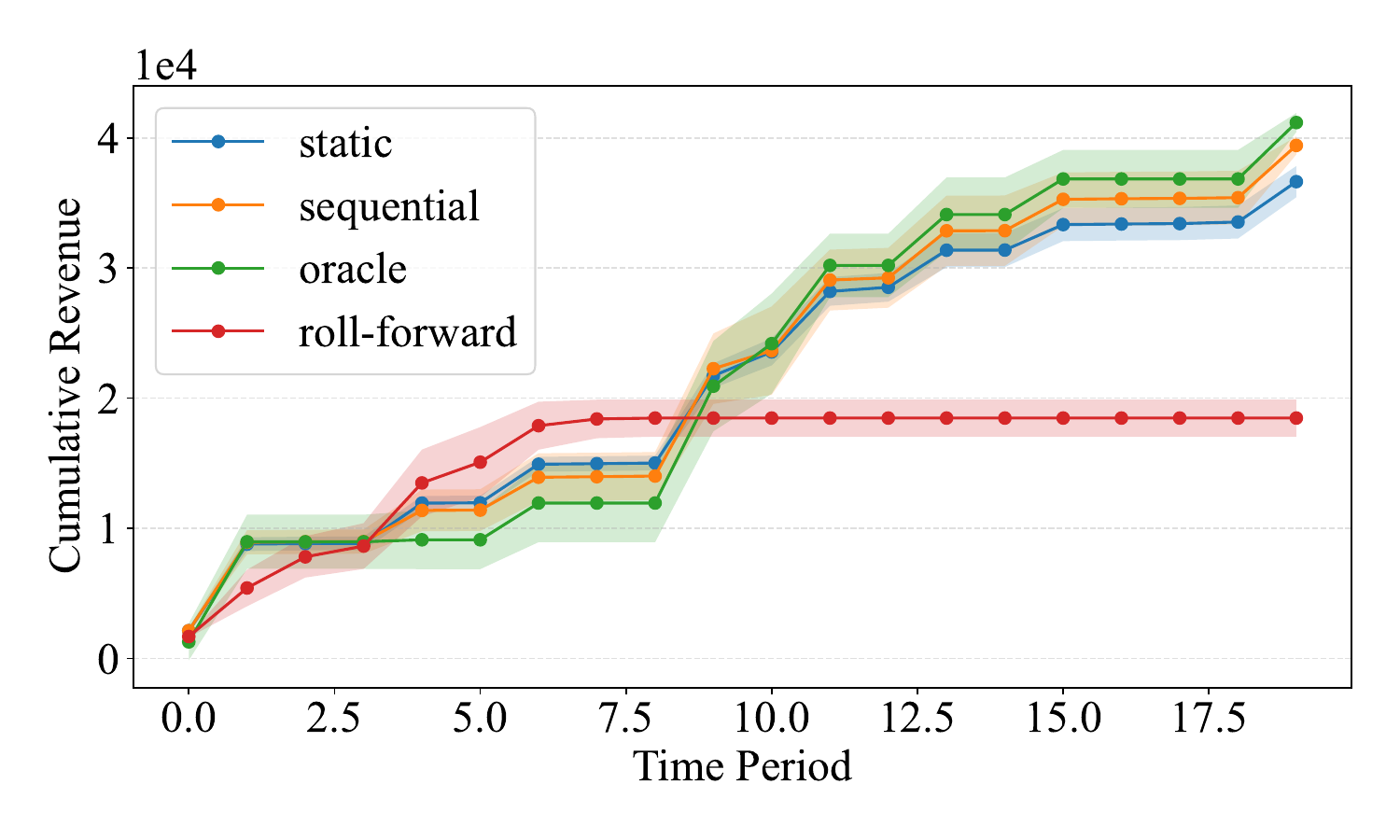}
    \caption{$T=20$.}
    \label{fig:stacked-a}
  \end{subfigure}
  \hfill
  \begin{subfigure}[t]{0.4\linewidth}
    \centering
    \includegraphics[width=\linewidth]{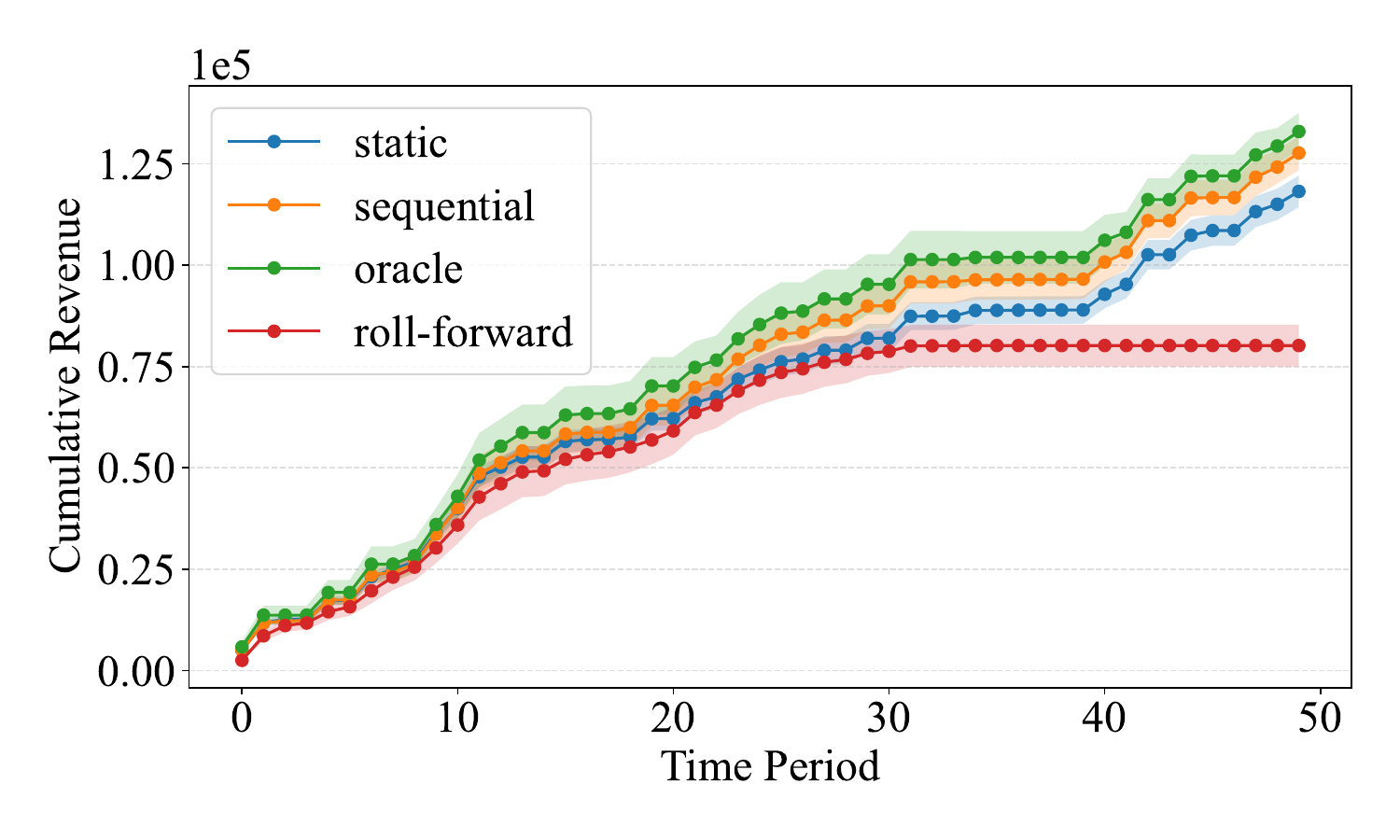}
    \caption{$T=50$.}
    \label{fig:stacked-b}
  \end{subfigure}

  \vspace{0.45em}

  \begin{subfigure}[t]{0.4\linewidth}
    \centering
    \includegraphics[width=\linewidth]{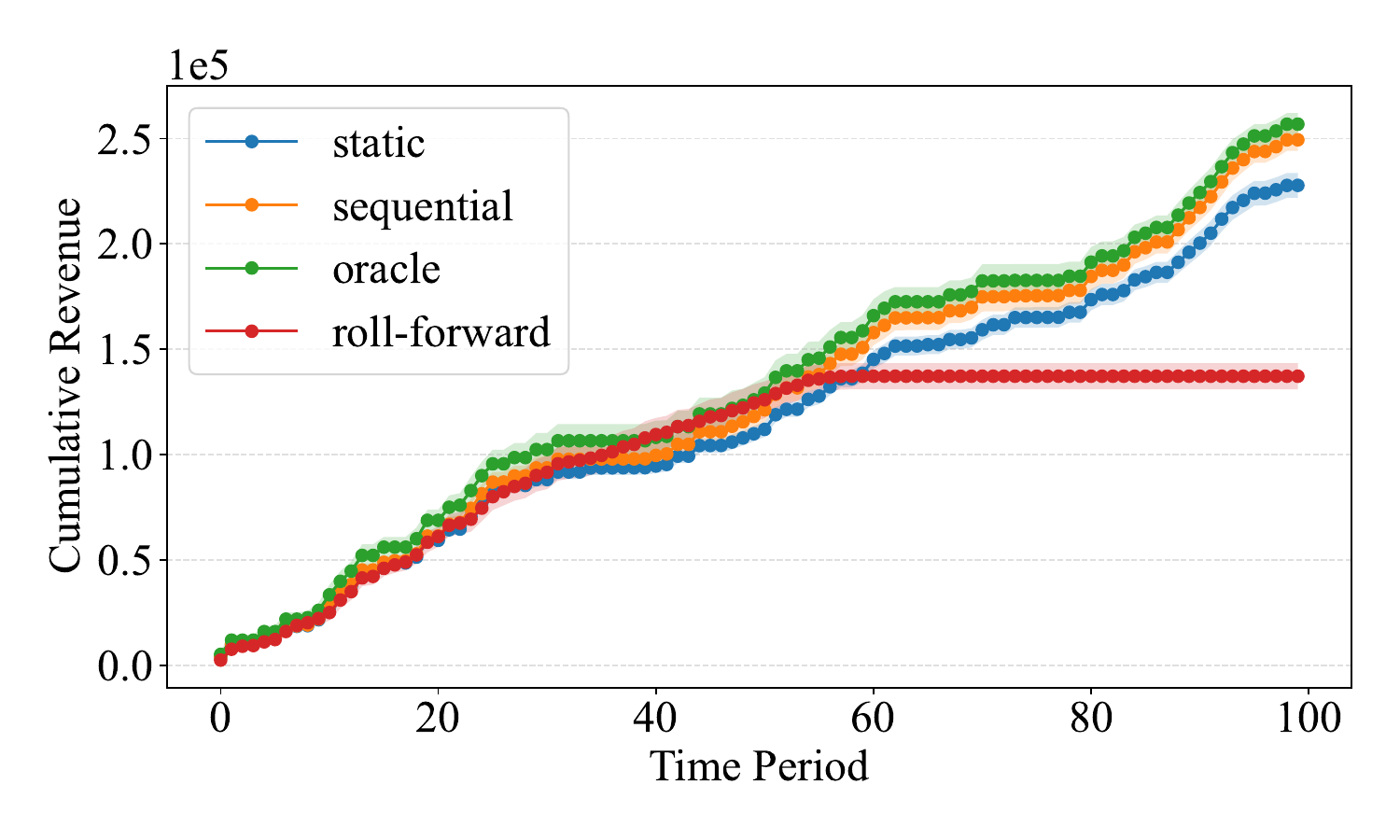}
    \caption{$T=100$.}
    \label{fig:stacked-c}
  \end{subfigure}
  \hfill
  \begin{subfigure}[t]{0.4\linewidth}
    \centering
    \includegraphics[width=\linewidth]{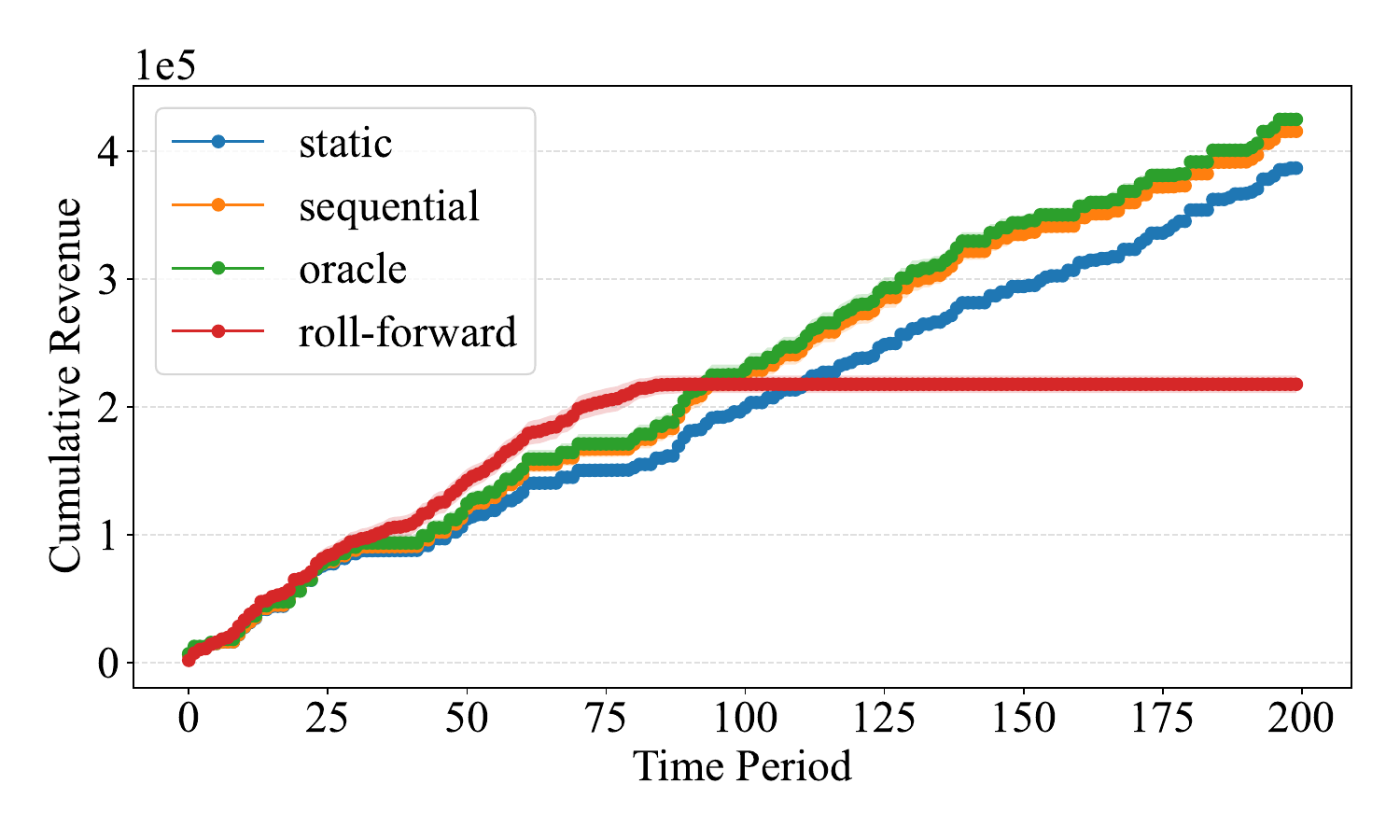}
    \caption{$T=200$.}
    \label{fig:stacked-d}
  \end{subfigure}

  \caption{The mean and 1-standard deviation interval for the cumulative revenue over $100$ trials for the various allocation methods and different horizons $T$.}
  \label{fig:aggreg_dep}
\end{figure*}

We include the aggregate results over $100$ trials for different values of $T$ in \Cref{tab:methods_by_T}. Across values of $T$ the \textit{sequential} allocation outperforms the \textit{static} one with respect to revenue. The \textit{roll-forward} allocation does not have a good performance. We show the mean and the 1-standard deviation interval for the cumulative revenue of the different allocation schemes in \Cref{fig:aggreg_dep} for the different values of $T$.

\newcommand{\tn}[2]{\shortstack{#1\\#2}}

\begin{table}[t]
  \centering
  \setlength{\tabcolsep}{7pt} % column padding; tweak as needed
  \renewcommand{\arraystretch}{1.6} % increase row height
  
  \begin{tabular}{lcccc}
    \toprule
    Method & $T=20$ & $T=50$ & $T=100$ & $T=200$ \\
    \midrule
    \textit{Roll-Forward}   & \tn{$18,473$}{($1,429$)} & \tn{$80,169$}{($5,154$)} & \tn{$137,179$}{($6,221$)} & \tn{$217,834$}{($6,804$)} \\
    \cmidrule(lr){1-5}
    \textit{Static}   & \tn{$36,644$}{($1,197$)} & \tn{$118,164$}{($3,892$)} & \tn{$227,780$}{($5,890$)} & \tn{$386,890$}{($5,615$)} \\
    \cmidrule(lr){1-5}
    \textit{Sequential} & \tn{$\mathbf{39,426}$}{($690$)} & \tn{$\mathbf{127,637}$}{($4,199$)} & \tn{$\mathbf{249,393}$}{($5,183$)} & \tn{$\mathbf{415,716}$}{($4,713$)} \\
    \cmidrule(lr){1-5}
    \textit{Oracle}   & \tn{$41,185$}{($734$)}  & \tn{$132,933$}{($4,566$)}  & \tn{$256,823$}{($5,306$)}  & \tn{$425,021$}{($5,011$)} \\
    \bottomrule
  \end{tabular}
  \caption{Results across time horizons $T$ for non-independent demands. Each cell shows two numbers (i.e., mean of achieved revenue on the first line and standard deviation of the achieved revenue in parentheses on the second line over $100$ trials). In bold, we emphasize the causal method that is best in terms of revenue, i.e., the \textit{sequential} allocation. We also observe that the \textit{sequential} allocation tends to have the smallest revenue variance. }
  \label{tab:methods_by_T}
\end{table}

\section{Conclusion}\label{sec:conclusion}

We propose a general sequential decision-making algorithm to solve the multi-period resource allocation problem under stochastic demands. Our approach is applicable under any model for stochastic demands and only requires knowledge of the conditional quantile functions. We show the superiority of our approach compared to the static solution using synthetic jointly log-normal data for the demands. The log-normal model is expressive and allows for easy calculation of the conditional quantile functions. Future work will focus on incorporating \Cref{alg:receding} in a closed loop, where the allocations are determined while the model for the distribution of the demands is learned. 
Furthermore, we plan to study in detail the importance of accurately modeling the distribution of the demands in applications.

%%%%%%%%%%%%%%%%%%%%%%%%%%%%%%%%%%%%%%%%%%%%%%%%%%%%%%%%%%%%%%%%%%%%%%%%%%%%%%%%

%===============================================================================

%===============================================================================

% The acknowledgments are automatically included only in the final and preprint versions of the paper.
\section*{Acknowledgments}

%Toyota Research Institute (TRI) provided funds to assist the authors with their research, but this article solely reflects the opinions and conclusions of its authors and not TRI or any other Toyota entity. 

%The NASA University Leadership Initiative (grant $\#$80NSSC20M0163) provided funds to assist the first author with their research, but this article solely reflects the opinions and conclusions of its authors and not any NASA entity. For the first author, 

The first author was partially funded through the Alexander S. Onassis Foundation Scholarship program. The authors would like to thank Logan Bell and Nikhil Devanathan for their feedback.

%%%%%%%%%%%%%%%%%%%%%%%%%%%%%%%%%%%%%%%%%%%%%%%%%%%%%%%%%%%%%%%%%%%%%%%%%%%%%%%%
%\bibliographystyle{IEEEtran}
%\bibliography{example}  % .bib
\renewcommand*{\bibfont}{\footnotesize}
\printbibliography

\end{document}